%% file: hicss.tex
\documentclass[10pt]{article}
\input{hicss-packages.tex}
\usepackage{hyperref}
\usepackage{multirow}

\usepackage{amsmath,amsfonts,amssymb,mathtools,nccmath}

\DeclareMathOperator*{\argmin}{arg\,min}

\newcommand{\di}{\textbf{i}}          

\title{AC-Network-Informed DC Optimal Power Flow for Electricity Markets}

\author{Gonzalo E. Constante \\
 Purdue University \\
 {geconsta@purdue.edu} \\ \And
 André H. Quisaguano \\
 Escuela Politécnica Nacional \\
 {andre.quisaguano@epn.edu.ec} \\ \And
  Antonio J. Conejo \\
 The Ohio State University \\
 {conejo.1@osu.edu} \\ \And
  Can Li\\
 Purdue University \\
 {canli@purdue.edu} \\
 }

\date{}

\begin{document}
\maketitle
\begin{abstract}

This paper presents a parametric quadratic approximation of the AC optimal power flow (AC-OPF) problem for time-sensitive and market-based applications. The parametric approximation preserves the physics-based but simple representation provided by the DC-OPF model and leverages market and physics information encoded in the data-driven demand-dependent parameters. To enable the deployment of the proposed model for real-time applications, we propose a supervised learning approach to predict near-optimal parameters, given a certain metric concerning the dispatch quantities and locational marginal prices (LMPs).  The training dataset is generated based on the solution of the accurate AC-OPF problem and a bilevel optimization problem, which calibrates parameters satisfying two market properties: cost recovery and revenue adequacy. We show the proposed approach's performance in various test systems in terms of cost and dispatch approximation errors, LMPs, market properties satisfaction, dispatch feasibility, and generalizability with respect to N-1 network topologies.

\end{abstract}

\subsubsection*{Keywords:}

Bilevel programming, electricity markets, machine learning, optimal power flow, parametric convex programming.

\section{Introduction}

Accurately modeling the governing physical laws of the flow of electricity over a network, also known as AC power flow (AC-PF) equations, is at the core of most of the decision-making processes of power systems, including control, operation, planning, and analysis. A cornerstone of power systems operations and electricity markets is the AC optimal power flow (AC-OPF), which aims to determine the most economical dispatch of generator output to meet the demand while satisfying the power flow equations and operational and physical constraints. However, the nonconvexity of the AC-PF equations raises critical computational challenges, especially for time-sensitive applications in large-scale systems \parencite{Bienstock2019}. Tractable surrogates of the AC-OPF problem have been proposed in the literature to circumvent these challenges, rendering easier-to-solve problems that can be solved in a timely manner \parencite{Molzahn2019}. 

The most common approximation of the AC-OPF is the DC optimal power flow (DC-OPF), whose linearity makes it suitable for a broad range of tasks, including market operations and power systems operations and long-term planning \parencite{Stott2009,Overbye2004}. However, the increasing electric demand from the transportation and data processing sectors is pushing power systems into conditions where the assumptions of the DC power flow approximation no longer hold. When these assumptions are not satisfied, the accuracy of DC power flow equations can be severely compromised, leading to inaccurate solutions that can result in an insecure and economically inefficient operation. As electricity demand continues to grow, the problem is exacerbated, further amplifying the impact of these violations.

One alternative to address this drawback is using machine learning (ML) techniques that leverage the solution of historical instances \parencite{VanHentenryck2021} to speed up or approximate the solution of the AC-OPF problem. ML-based approaches for power systems and markets operations include finding warm starting points \parencite{Baker2019}, learning active sets \parencite{Misra2022,Ferrando2024}, and end-to-end learning frameworks \parencite{Zamzam2020,Chen2022}. ML models can be enhanced by targeting or enforcing the satisfaction of physical laws or engineering constraints \parencite{Ng2018} for power systems operational tools using soft and hard approaches. Soft approaches, which augment the loss function with penalty terms to enforce constraints on the primal/dual variables, include models informed by physics \parencite{Fioretto2020,Lei2021,Nellikkath2022}, 
topology \parencite{Falconer2023}, and sensitivities \parencite{Singh2022}. On the other hand, hard approaches aim to strictly satisfy the constraint set or a subset of it during the inference stage \parencite{Pan2023,Li2023}.

Another alternative to represent AC power flow equations is to use tractable parametric relaxations/approximations of such equations, whose parameters are fine-tuned to minimize the error with respect to a given metric, e.g., dispatch accuracy. 
Recent works by \cite{Taheri2023,Taheri2024} propose an unconstrained optimization model to find a set of multiplicative and additive parameters, which are optimized for a range of operating conditions, to improve the accuracy of the DC power flow model with respect to its AC counterpart. Such framework preserves the approximation and simplicity provided by the DC power flow equations while refining it with fine-tuned parameters. In another context, \cite{Gupta2024} propose a bilevel model to learn parameters for a tractable surrogate of complex optimization problems, which is solved using an exact penalty reformulation. In the above works, once the parameters are optimized for a given dataset representing different operating conditions, the optimized parameters remain fixed during the inference (online) stage for all operating conditions.

In this paper, we propose a parametric optimization model for electricity market applications that relies on the physics-based simple model provided by the DC-OPF model and refines it using data-driven information encoded in a set of tunable parameters. Such parameters are determined to reduce the dispatch approximation error of the parametric model with respect to the AC-OPF model while targeting desirable market design properties. Additionally, the proposed parametric model is interpretable since the tunable parameters correspond to scaling factors of the nodal demands aimed at correcting the errors introduced due to the assumptions of the DC approximation. 

The proposed approach entails an offline and an online stage. In the offline (training) stage, we propose a supervised learning approach to learn the mapping between active and reactive demand profiles and the tunable parameters. The optimal tunable parameters are determined based on a bilevel model. The upper-level problem of the bilevel model minimizes the dispatch error between the AC-OPF problem and its parametric DC counterpart while satisfying cost recovery and revenue adequacy. The lower-level problem of the bilevel model is the proposed parametric DC-OPF problem. In the online stage, we predict (near-)optimal tunable parameters and solve the parametric DC-OPF problem.

The remainder of the paper is organized as follows. \hyperref[sec:background]{Section~2} introduces the AC- and DC-OPF models, market properties, and the proposed parametric model. \hyperref[sec:methodology]{Section~3} discusses the proposed learning-based methodology. \hyperref[sec:experiments]{Section~4} presents numerical experiments from representative test cases assessing the performance of the proposed model. In \hyperref[sec:conclusion]{Section~5}, we present the concluding remarks and future work.

\section{Parametric DC-OPF Problem} \label{sec:background}

This section presents the AC- and DC-OPF problem formulations and market properties in the context of the common marginal pricing scheme. Then, we present the proposed parametric model and a bilevel optimization model to determine the optimal values of the tunable parameters.

We first establish our notation. Let $\mathcal{N}$ denote the set of nodes, and $\mathcal{E} \subseteq \mathcal{N}\times \mathcal{N}$ denote the set of branches. Every node $i$ is associated with a voltage phasor with magnitude $v_i$, a phase angle $\theta_i$, and a shunt admittance $Y_{i}^{\rm sh} = G_{i}^{\rm sh} + \di B_{i}^{\rm sh}$ where $\di = \sqrt{-1}$. We denote the active and reactive power demand at node $i$ as $P_i^{\rm d}$ and $Q_i^{\rm d}$, respectively. Each branch is indexed by $(i,j)$ and characterized by a series admittance $Y_{ij} = G_{ij} + \di B_{ij}$, a shunt susceptance $B_{ij}^{\rm sh}$, a complex power flow $s_{ij} = p_{ij} + \di q_{ij}$, and a capacity limit $\overline{S}_{ij}$. Let $\mathcal{G}$ be the set of generators indexed by $k$ whose active and reactive power production are denoted by $p^{\rm{g}}_k$ and $q^{\rm{g}}_k$, respectively. Let $\mathcal{G}_i$ denote the set of generators connected to node $i$.
Lowercase and uppercase letters represent optimization variables and parameters, respectively, whereas boldface letters denote matrices and vectors.

\subsection{AC Optimal Power Flow}

The AC-OPF problem can be formulated as follows:

\begin{subequations} \label{eq:1}
\begin{align}
\min_{\Xi^{\rm{AC}}} \quad & \sum_{k \in \mathcal{G}} \left(C^{(2)}_k \left(p_k^{\rm g}\right)^2 + C^{(1)}_k p_k^{\rm g} \right) \label{eq:1a} \\
\textrm{s.t.\quad}   & (\forall i \in \mathcal{N}, \forall (i,j) \in \mathcal{E}, \forall k \in \mathcal{G}): \nonumber \\
& \sum_{k\in\mathcal{G}_i}p_k^{\rm g} -P_i^{\rm d} - G_{i}^{\rm sh}c_{ii} = \smashoperator{\sum_{(i,j)\in\mathcal{E}}}p_{ij},  \label{eq:1b} \\
& \sum_{k\in\mathcal{G}_i}q_k^{\rm g} -Q_i^{\rm d} + B_{i}^{\rm sh}c_{ii} = \smashoperator{\sum_{(i,j)\in\mathcal{E}}}q_{ij}, \label{eq:1c} \\
& p_{ij} = G_{ij}c_{ii} - G_{ij}c_{ij} + B_{ij}s_{ij}, \label{eq:1d} \\
& \begin{multlined}
    q_{ij} = -\left(B_{ij} + B_{ij}^{\rm sh}\right)c_{ii} + G_{ij}s_{ij} + \\  B_{ij}c_{ij},
\end{multlined}  \label{eq:1e} \\
& \left(p_{ij}\right)^2 + \left(q_{ij}\right)^2 \le \overline{S}_{ij}^2,  \label{eq:1f} \\
& \underline{V}_i^2 \le c_{ii} \le \overline{V}_i^2,     \label{eq:1g} \\
& \underline{P}_k \le p_k^{\rm g} \le \overline{P}_k,\label{eq:1h} \\
& \underline{Q}_k \le q_k^{\rm g} \le \overline{Q}_k, \label{eq:1i} \\
& c_{ii} = v_i^2, \label{eq:1j} \\
& c_{ij} = v_i v_j\cos{(\theta_i -\theta_j)}, \label{eq:1k} \\
\phantom{\textrm{s.t.\quad}}\;& s_{ij} = -v_i v_j\sin{(\theta_i -\theta_j)}, \label{eq:1l} \qquad \qquad \quad
\end{align}
\end{subequations}
where the optimization variables are elements of the set
\[
\Xi^{\rm{AC}} = \{p_k^{\rm g},q_k^{\rm g},p_{ij},q_{ij},c_{ii},c_{ij},s_{ij},v_i,\theta_i\},
\]
and $C^{(2)}_k$ and $C^{(1)}_k$ denote the quadratic and linear cost coefficients of generator $k$. 

The objective function \eqref{eq:1a} is the total production cost. Constraints \eqref{eq:1b} and \eqref{eq:1c} enforce the active and reactive power balance, respectively. Constraints \eqref{eq:1d}-\eqref{eq:1e} define the active and reactive power flows throughout the power grid. Constraints \eqref{eq:1f} enforce the thermal capacity of the branches. Constraints \eqref{eq:1g} bound the squared voltage magnitudes. Constraints \eqref{eq:1h}-\eqref{eq:1i} bound the generating units' active and reactive power output, respectively. Constraints \eqref{eq:1j}-\eqref{eq:1l} are the nonconvex relationships between the voltage magnitudes and phase angles and the lifted variables, $c_{ii}$, $c_{ij}$ and $s_{ij}$. 

\subsection{DC Optimal Power Flow}

In most operating conditions, bulk power systems satisfy the following criteria: (i) angle differences between two connected nodes are small, (ii) voltage magnitudes are near their nominal values, and (iii) branches have a high reactance/resistance ratio. Using the criteria above, the AC-OPF problem can be approximated by the DC-OPF problem, which can be formulated as follows:

\begin{subequations} \label{eq:2}
\begin{align}
\min_{\Xi^{\rm{DC}}} \quad & \sum_{k \in \mathcal{G}} \left(C^{(2)}_k \left(p_k^{\rm g}\right)^2 + C^{(1)}_k p_k^{\rm g} \right) \label{eq:2a} \\
\textrm{s.t.} \quad  & (\forall i \in \mathcal{N}, \forall (i,j) \in \mathcal{E}, \forall k \in \mathcal{G}) \nonumber \\
& \sum_{k\in\mathcal{G}_i}p_k^{\rm g} -P_i^{\rm d} - G_{i}^{\rm sh} = \smashoperator{\sum_{(i,j)\in\mathcal{E}}}p_{ij}, \quad (\lambda_i),  \label{eq:2b} \\
& p_{ij} = B_{ij}\left(\theta_i - \theta_j \right), \label{eq:2c} \\
& -\overline{S}_{ij} \le p_{ij} \le \overline{S}_{ij},  \label{eq:2d} \\
& \underline{P}_k \le p_k^{\rm g} \le \overline{P}_k,\label{eq:2e}
\end{align}
\end{subequations}

\noindent
where the optimization variables are elements of the set $\Xi^{\rm{DC}}=\{p_k^{\rm g},p_{ij},\theta_i\}$. Constraints \eqref{eq:2b} and \eqref{eq:2c} correspond to the linear approximation of the power balance and power flow equations, respectively. Constraints \eqref{eq:2d} and \eqref{eq:2e} bound the power flows and generating units' active power generation, respectively. The dual variables, $\lambda_i$ of the power balance constraints \eqref{eq:2b}, are called locational marginal prices (LMP) and are used for market pricing purposes.

\subsection{Energy Marginal Pricing and Market Design Properties}

In electricity markets, LMPs reflect the cost of supplying the next marginal increment of load at different locations considering generation costs and the effect of power losses and constraints of the transmission network \parencite{Tanaka2022}. LMPs along with the generation dispatch enable economic transactions and should ensure the satisfaction of key principles, i.e., cost recovery and revenue adequacy, to achieve an economically efficient power system operation.

\subsubsection{Cost Recovery}

A market complies with the cost recovery property when the profit of each producer is nonnegative. From the operation perspective, the payments made to each producer based on the corresponding LMP must cover its operating costs. The cost recovery property can be formulated as follows:
\begin{equation}
\lambda_{i(k)} p_k^{\rm g} \ge C^{(2)}_k \left(p_k^{\rm g}\right)^2 + C^{(1)}_k p_k^{\rm g},\; \forall k \in \mathcal{G},
\end{equation}
where $\lambda_{i(k)}$ denotes the LMP of the node of generator $k$. In particular, the above condition is satisfied with equality for the generators whose production levels do not reach their upper/lower bounds. 

\subsubsection{Revenue Adequacy}

A market is said to be revenue adequate if the market operator collects from the consumers at least the amount needed to pay the producers. In a marginal pricing setting, consumers and producers are charged and paid, respectively, based on their corresponding LMPs. Hence, the revenue adequacy condition can be formulated as follows:
\begin{equation}
\sum_{i\in \mathcal{N}}\lambda_i P_i^{\rm d} \ge \sum_{k\in \mathcal{G}}\lambda_{i(k)} p_k^{\rm g}
\end{equation}

\subsection{Proposed Model Formulation}

The increased electric demand of data centers and electrification of industrial processes and transportation are driving the operation of power systems to operating points where DC approximation conditions are not generally satisfied. Despite being the cornerstone of power system operations and electricity markets, the errors in the generation dispatch of the DC-OPF problem, especially in such operating conditions can result in an insecure and economically inefficient operation. To address this issue, we propose a parametric approximation of the AC-OPF problem, leveraging the physics-based power flow simple approximation of the DC-OPF problem with data-driven methods to reduce dispatch errors. The proposed pDC-OPF problem can be formulated as follows:
\begin{subequations} \label{eq:5}
\begin{align}
\min_{\Xi^{\rm{DC}}} \quad & \sum_{k \in \mathcal{G}} \left(C^{(2)}_k \left(p_k^{\rm g}\right)^2 + C^{(1)}_k p_k^{\rm g} \right) \nonumber \\
\textrm{s.t.} \quad  & \sum_{k\in\mathcal{G}_i}p_k^{\rm g} - P_i^{\rm d} \beta_i\left(\boldsymbol{P}^{\rm d},\boldsymbol{Q}^{\rm d}\right) -  G_{i}^{\rm sh} = \nonumber \\
& \quad \smashoperator{\sum_{(i,j)\in\mathcal{E}}}p_{ij}, \quad  \forall i \in \mathcal{N}, \quad (\lambda_i),   \label{eq:5b} \\
& \eqref{eq:2b} - \eqref{eq:2d},
\end{align}
\end{subequations}
where $\beta_i(\boldsymbol{P}^{\rm d},\boldsymbol{Q}^{\rm d})$ denotes the scaling parameter at node $i$, which encodes information about active and reactive power demand at every node, and $\lambda_i$  denotes the corresponding LMP. Note that $\beta_i$ scales the demand at node $i$ to reduce the errors induced by the DC-OPF approximation.

The solution of the proposed pDC-OPF model provides a generation dispatch that can be used as a surrogate of the AC-OPF model for power systems operations. To enable the pDC-OPF model as a possible market tool, we should target that the scaling parameter $\beta_i(\boldsymbol{P}^{\rm d},\boldsymbol{Q}^{\rm d})$ renders a dispatch and LMPs that satisfy the cost recovery of all generating units and revenue adequacy properties. It is important to note that the latter property needs to be satisfied considering the original (unscaled) demand. 

\subsection{Determining Optimal Tunable Parameters}

For a given demand pair $(\boldsymbol{P}^{\rm d},\boldsymbol{Q}^{\rm d})$, the optimal tunable parameters, satisfying the cost recovery and revenue adequacy properties, can be determined by solving the following bilevel optimization problem:
\begin{subequations} \label{eq:7}
\begin{align}
\min_{\beta_i\ge 0} \quad & \frac{\gamma^{\rm p}}{\lvert \mathcal{G} \rvert}\sum_{k\in\mathcal{G}} \left(p^{\rm g}_k - p^{\rm g,*}_k\right)^2 + \nonumber \\
& \frac{\gamma^{\lambda}}{\lvert \mathcal{N} \rvert}\sum_{i\in\mathcal{N}} \left(\lambda_i - \lambda_i^*\right)^2 +  \frac{\gamma^{\beta}}{\lvert \mathcal{N} \rvert}\sum_{i\in\mathcal{N}} \left(\beta_i\right)^2\label{eq:7a} \\
\textrm{s.t.} \quad  & \lambda_{i(k)} p_k^{\rm g} \ge c^{(2)}_k \left(p_k^{\rm g}\right)^2 + c^{(1)}_k p_k^{\rm g},\; \forall k \in \mathcal{G}, \label{eq:7b} \\
& \sum_{i\in \mathcal{N}}\lambda_i P_i^{\rm d} \ge \sum_{k\in \mathcal{G}}\lambda_{i(k)} p_k^{\rm g}, \label{eq:7c} \\
& \min_{\Xi^{\rm{DC}}} \quad \sum_{k \in \mathcal{G}} \left(c^{(2)}_k \left(p_k^{\rm g}\right)^2 + c^{(1)}_k p_k^{\rm g} \right) \label{eq:7d} \\
& \;\,\textrm{s.t.} \quad (\forall i \in \mathcal{N}, \forall (i,j) \in \mathcal{E}, \forall k \in \mathcal{G}) \nonumber \\
& \phantom{\;\,\textrm{s.t.}} \quad \sum_{k\in\mathcal{G}_i}p_k^{\rm g} - P_i^{\rm d} \beta_i -  G_{i}^{\rm sh} = \nonumber \\
& \phantom{\;\,\textrm{s.t.}} \quad \quad \smashoperator{\sum_{(i,j)\in\mathcal{E}}}p_{ij}, \quad (\lambda_i),   \label{eq:7e} \\
& \phantom{\;\,\textrm{s.t.}} \quad \eqref{eq:2b}-\eqref{eq:2d}, \nonumber
\end{align}
\end{subequations}
where $\gamma^{\rm p}$, $\gamma^{\lambda}$, and $\gamma^{\beta}$ are weighting coefficients, $p^{\rm g,*}_k$ denotes the optimal dispatch of unit $k$ for the AC-OPF problem, and $\lambda_i^*$ denotes the dual variable of the $i$-th active power balance constraint \eqref{eq:1b} at the optimal solution. 

The objective function \eqref{eq:7a} of the upper-level problem aims at minimizing a trade-off of (i) the mean squared error between the dispatch and LMPs of the pDC-OPF problem and the corresponding ones of the AC-OPF problem and (ii) an $
\ell_2$-penalty to $\boldsymbol{\beta}$. Such objective function will promote a set of small parameters due to the shrinkage effect of the $\ell_2$-norm regularization, with reduced dispatch and LMPs errors.  Including such regularization promotes smoothness and reduces the sparsity and discontinuity of the optimal tunable parameters, which improves the training process described in Section \ref{sec:methodology}.
Upper-level constraints \eqref{eq:7b} and \eqref{eq:7c} enforce the cost recovery and revenue adequacy properties, respectively, of the solution of the pDC-OPF problem. The lower-level problem, \eqref{eq:7d}, \eqref{eq:7e} and \eqref{eq:2b}-\eqref{eq:2d}, models the pDC-OPF problem.

Problem \eqref{eq:7} can be solved by replacing the lower-level problem by (i) its KKT-optimality conditions, which can then be linearized using a Fortuny-Amat \& McCarl linearization \parencite{Gabriel2013}, or by (ii) its primal and dual constraints and enforcing the strong duality equality, which involves products of $\beta_i$ and the dual variables $\lambda_i$ \parencite{Ruiz2009}. Both alternatives result in computationally challenging problems, which may be impractical for real-time applications. 

In the following section, we present a supervised learning approach to circumvent this key issue by learning, in an offline stage, the mapping between the demand vector pair $(\boldsymbol{P}^{\rm d},\boldsymbol{Q}^{\rm d})$ and their optimal tunable parameters $\boldsymbol{\beta}^*$, and predicting, at inference time, near-optimal tunable parameters $\hat{\boldsymbol{\beta}}$, which are then used in the proposed pDC-OPF model.

\section{Proposed Approach} \label{sec:methodology}

\begin{figure}[t]
    \centering
    \includegraphics[width=\linewidth]{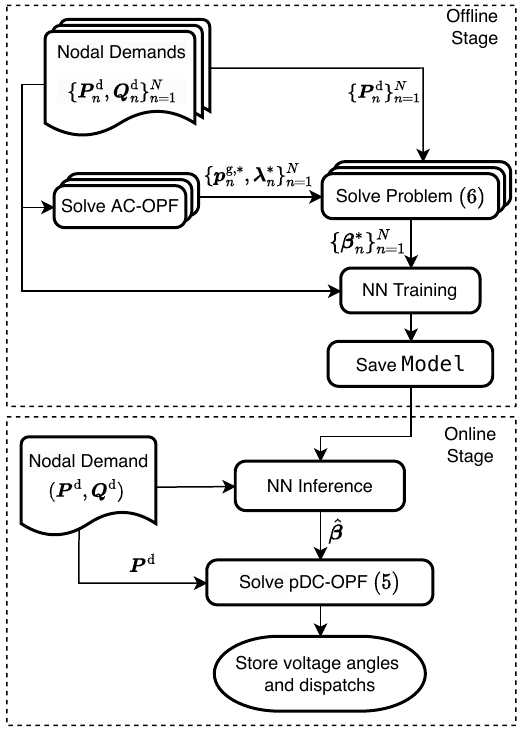}
    \caption{Proposed approach flowchart}
    \label{fig:flowchart}
\end{figure}

The proposed approach aims to enable the deployment of the pDC-OPF problem as a tool for time-sensitive market operations. This approach, illustrated in Figure~\ref{fig:flowchart}, involves an offline (training) stage and an online (execution) stage. The offline stage is focused on training a neural network (NN) in a supervised learning fashion to predict (near-)optimal tunable parameters for a given demand. In the online stage, the parameters are predicted and the pDC-OPF problem \eqref{eq:5} is solved.

\subsection{Offline Stage}

In the offline stage, we train a neural network using a supervised learning approach to learn the mapping between active and reactive nodal demands and optimal scaling vector, $\boldsymbol{\beta}^\ast$. We first generate $N$ samples of active and reactive power demand vectors and the solution of the AC-OPF problem, i.e., a tuple $\{\boldsymbol{P}^{\rm d}_n,\boldsymbol{Q}^{\rm d}_n,\boldsymbol{p}^{\rm g,*}_n,\boldsymbol{\lambda}^*_n\}_{n=1}^N$, where $\boldsymbol{P}^{\rm d}_n$ and $\boldsymbol{Q}^{\rm d}_n$ are the $n$-th active and reactive power demand vectors, and $\boldsymbol{p}^{\rm g,*}_n$ and $\boldsymbol{\lambda}^*_n$ the dispatch and LMP vectors at the solution of the corresponding AC-OPF problem. Then, for each sample, we solve the bilevel optimization problem \eqref{eq:7} to compute the optimal scaling vector $\boldsymbol{\beta}_n^*$. 

The NN training dataset contains the above $N$ samples of active and reactive power demand vectors and the corresponding optimal scaling vectors. Then, we train a NN model to learn the mapping:
\begin{align}
\hat{\boldsymbol{\beta}}_n = \mathrm{NN} (\boldsymbol{x}_n;\boldsymbol{\Theta}), 
\end{align}
where $\hat{\boldsymbol{\beta}}_n \approx \boldsymbol{\beta}_n^*$ is the predicted scaling vector for the $n$-th sample, $\Theta$ denotes the parameters (weights and biases) of the NN, and $\boldsymbol{x}_n = [\boldsymbol{P}^{\rm d}_n,\boldsymbol{Q}^{\rm d}_n]$ denotes the $n$-th input sample. We consider the following loss function to train the NN:
\begin{align}
    \mathcal{L}\left(\boldsymbol{\Theta}\right) = &
    \frac{1}{N}\sum_{n=1}^{N} \lVert \mathrm{NN} (\boldsymbol{\Theta}; \boldsymbol{x}_n) \odot \boldsymbol{P}^{\textrm d}_n - \boldsymbol{\beta}_n^* \odot \boldsymbol{P}^{\textrm d}_n \rVert^2 + \nonumber \\
    & \frac{\rho}{N} \sum_{n=1}^{N} \lVert \mathrm{NN} (\boldsymbol{\Theta}; \boldsymbol{x}_n)^\top \boldsymbol{P}_n^{\rm d} - \left(\boldsymbol{\beta}_n^*\right)^\top \boldsymbol{P}^{\textrm d}_n \rVert^2, \label{eq:loss}
\end{align}
where $\odot$ denotes the Hadamard (element-wise) product and $\rho$ is the regularization parameter. The first term of the loss function \eqref{eq:loss} corresponds to the mean squared error between the optimal scaled demand and the predicted one whereas the second term penalizes the difference between the total scaled demand and the total predicted one. Once the NN is trained, we save its parameters such that the NN can be deployed to predict the scaling vector for real-time applications.

\subsection{Online Stage}

In the online stage, we first predict a near-optimal scaling vector $\hat{\boldsymbol{\beta}}$ for a given pair of active and reactive power demands, $(\boldsymbol{P}^{\rm d},\boldsymbol{Q}^{\rm d})$ using the NN. Then, the pDC-OPF model \eqref{eq:5} is solved for the predicted scaling vector.

The proposed model's performance depends on the accuracy of the trained NN. However, it is important to note that since the proposed parametric formulation is agnostic to the prediction model, one can use any machine learning architecture, loss function, or training algorithm to learn the mapping between nodal demands and optimal scaling vector.

\section{Numerical Experiments} \label{sec:experiments}

\subsection{Experimental Setup}
In this section, we provide numerical experiments to illustrate the performance of the proposed pDC-OPF model. The numerical experiments have been implemented on a Lenovo Y520 laptop with an Intel Core i7-7700HQ processor clocking at 2.80GHz and 16 GB of RAM. Our numerical experiments are conducted using three IEEE test systems of 30, 57, and 118 buses, whose nominal parameters can be retrieved from the Power Grid Library (PGLib-OPF) v23.07  \parencite{Babaeinejadsarookolaee2019}. The code is available in \url{https://github.com/li-group/AC-Network-Informed-DC-Optimal-Power-Flow}.

In the offline stage, we create the dataset of 7,500 uniformly sampled demands for the 30- and 57-bus systems and 18,500 samples for the 118-bus system with a $80\%/20\%$ split for training/testing for each system. The active and reactive power demands are sampled within $[70\%, 130\%]$ and $[85\%, 100\%]$, respectively, of the original load available in the PGLib-OPF instances. The solutions of the AC-OPF and DC-OPF problems are generated in MATLAB 2023a using Matpower \parencite{Zimmerman2011}, whereas the bilevel optimization problem and the pDC-OPF model have been implemented under JuMP 1.10 with Knitro as the local optimization solver. 

The BilevelJuMP package \parencite{DiasGarcia2023} has been used to implement and recast the bilevel problem based on its primal and dual constraints and the strong duality equality. To focus on reducing the dispatch errors instead of prices accuracy, we set the weight coefficients of the upper-level objective function to $\gamma^{\lambda} = 0$ and $\gamma^{\beta} = 1$. Also, we fix the dispatch of the lower-level problem to the dispatch of the AC-OPF problem to speed up the solution of Problem \eqref{eq:7}. 

The NNs were trained in Python 3.9 using the PyTorch environment, with the Adam optimizer set to a learning rate of $1\mathrm{e}{-5}$ and a weight decay of $1\mathrm{e}{-4}$. For each test system, we train a NN with two hidden layers with rectified linear units (ReLU) activations, respectively, and a linear output layer. For the 30- and 57-bus systems, the hidden layers have 512 and 256 neurons, whereas for the 118-bus system, they have 1024 and 512 neurons. The training process spanned 500 epochs with a batch size of 64 samples.

We generate 1,000 additional demand samples, which are not in the training/testing datasets, to validate the proposed model's out-of-sample performance in terms of cost error approximation, dispatch approximation and feasibility, and satisfaction of market design properties. Additionally, we assess the ability of the proposed model to generalize with respect to N-1 network topologies, i.e., its ability to perform well with respect to N-1 topologies that are not represented in the dataset of the training stage.

\subsection{Performance Analysis}

We assess the proposed model's performance using several metrics to measure the cost and dispatch approximation errors and how close the dispatch is to satisfying the constraints of the AC-OPF problem.
\paragraph{Total cost} We define the following metric to quantify the cost approximation error with respect to the one of the AC-OPF problem:
\begin{equation}
\text{Cost error}\,(\%) = 100 \cdot \frac{\lvert {\rm Cost}^{\rm AC} - {\rm Cost}^{\rm X} \rvert}{{\rm Cost}^{\rm AC}} \label{eq:cost_error}
\end{equation}
where ${\rm X} = \{\text{DC-OPF, pDC-OPF}\}$. Table~\ref{tab:cost} presents the maximum, mean, and minimum errors of the pDC-OPF and DC-OPF model. The pDC-OPF model shows a better performance in all the statistics for the three studied systems; however, the error of the DC-OPF model has a smaller range of variation. We also note that the maximum errors of the proposed model are better than the minimum ones of the DC-OPF model, i.e., our model showed smaller errors for all the generated samples of the three systems.

\begin{table}[t]
\centering
\caption{Cost approximation errors}
\label{tab:cost}
\begin{tabular}{c|l|lll}
\hline
\multirow{2}{*}{\textbf{Model}}   & \multirow{2}{*}{\textbf{Stat.}} &   \multicolumn{3}{|c}{\textbf{Cost error (\%)}} \\
                         &   & \bf  30-bus   & \bf  57-bus   & \bf 
 118-bus   \\
\hline
\bf  \multirow{3}{*}{pDC-OPF} & min   & 0.0004   & 0.0004   & 0.0004    \\
                         & mean  & 0.2607   & 0.2441   & 0.4205    \\
                         & max   & 0.8612   & 0.9998   & 1.9695    \\
\hline
\bf \multirow{3}{*}{DC-OPF}  & min   & 8.0823   & 5.7582   & 3.5567    \\
                         & mean  & 8.6006   & 7.2994   & 4.1807    \\
                         & max   & 9.5172   & 8.4568   & 4.8353    \\
\hline
\end{tabular}
\end{table}

\begin{figure}[t]
    \centering
    \includegraphics[width=0.825\linewidth]{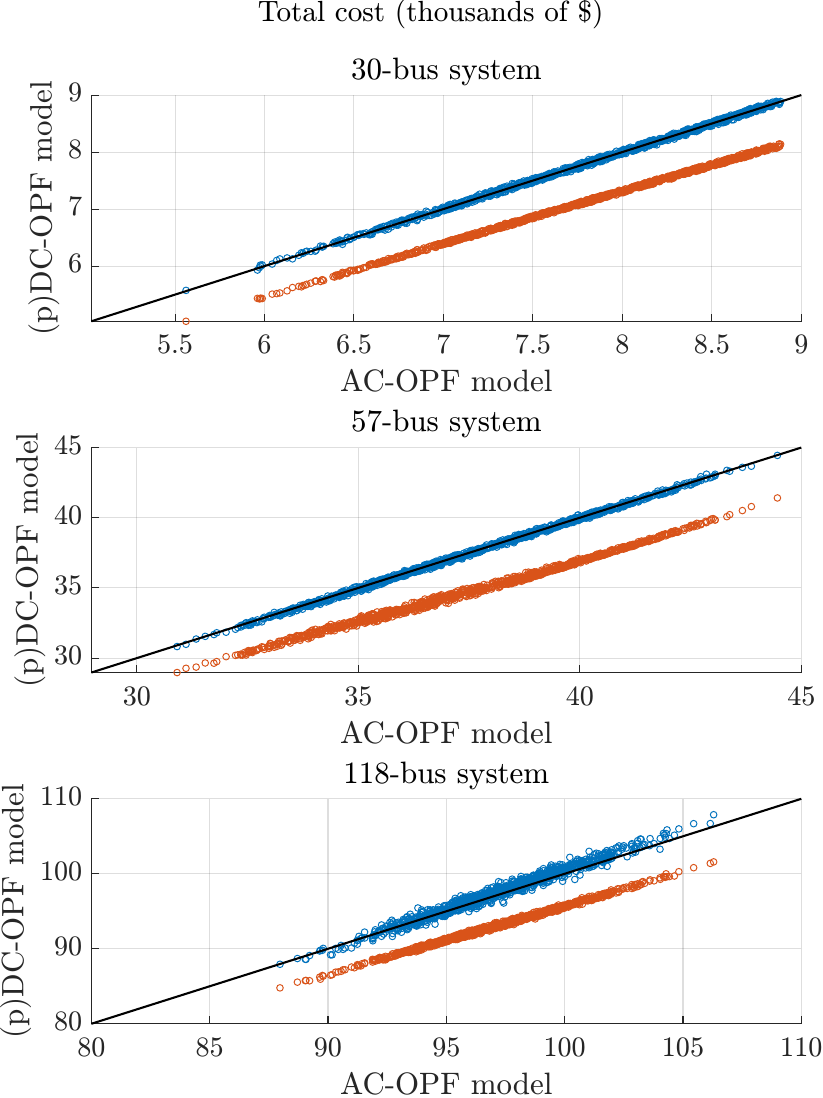}
    \caption{Cost correlation with respect to AC-OPF model. DC-OPF (red). pDC-OPF (blue).}
    \label{fig:CostComparison}
\end{figure}

Figure~\ref{fig:CostComparison} illustrates the cost correlation of the (p)DC-OPF models with respect to the AC-OPF model. The proposed model clearly outperforms the DC-OPF model for all the samples for the three systems. We note that the cost offset of the DC-OPF model is due to neglecting system losses, which are implicitly considered in the parameters of the proposed model.

\paragraph{Generation dispatch}

We also study the performance of the proposed model in terms of (i) the dispatch optimality error, i.e., the distance with respect to the optimal dispatch of the AC-OPF problem, and (ii) the distance to feasibility, i.e., the error of the dispatch with respect to its closest projection onto the nonconvex set.

\begin{figure}[t]
    \centering
    \includegraphics[width=\linewidth]{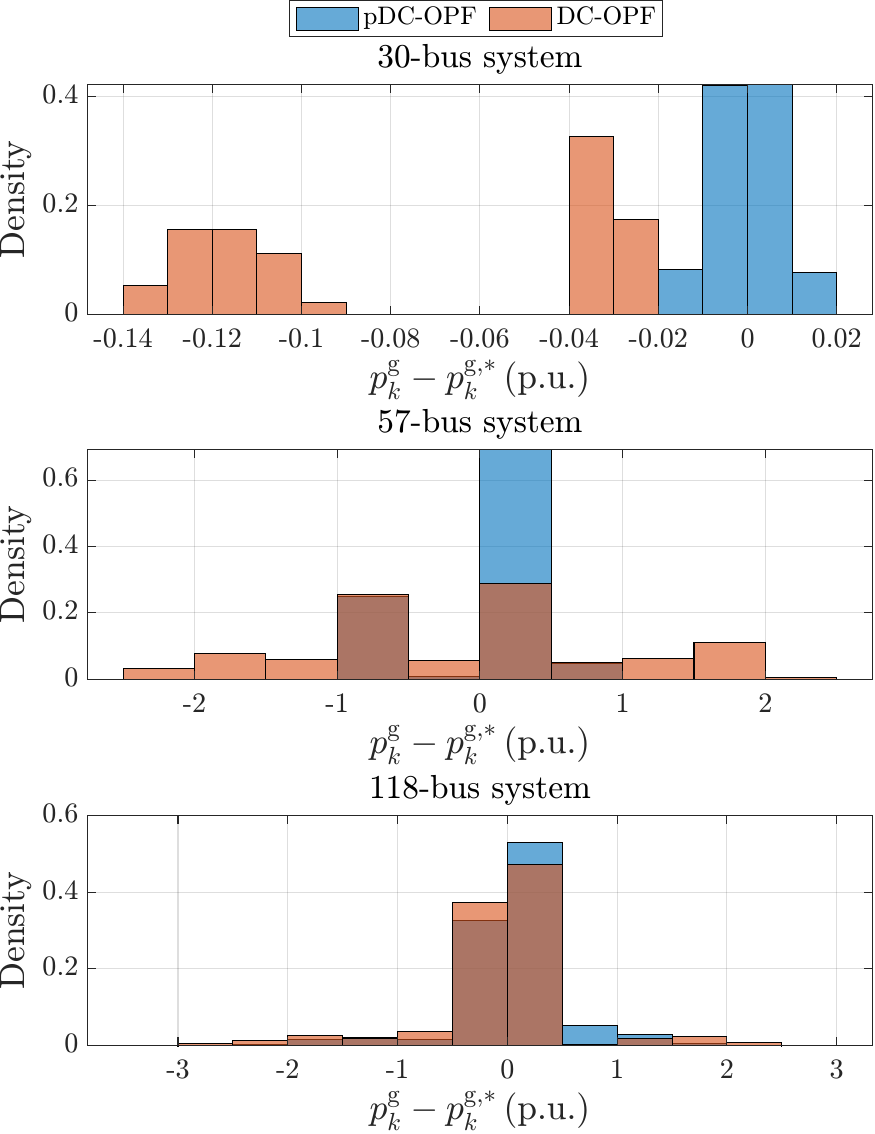}
    \caption{Dispatch optimality error.}
    \label{fig:dispatch_error}
\end{figure}

Figure~\ref{fig:dispatch_error} illustrates the dispatch optimality error. For the three test systems, we note that the proposed model improves the mean and worst-case dispatch errors with respect to the original DC model.

\begin{figure}[t]
    \centering
    \includegraphics[width=\linewidth]{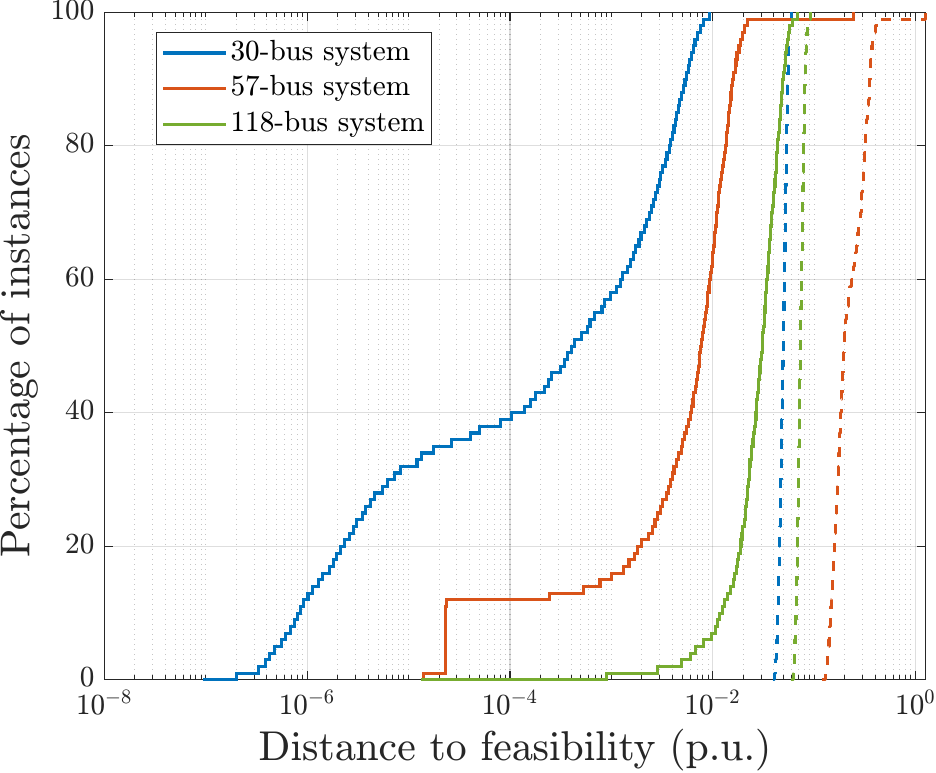}
    \caption{Distance to feasibility. pDC-OPF model (solid). DC-OPF model (dashed).}
    \label{fig:Feasibility}
\end{figure}

To measure how close the dispatch of the generators of the (p)DC models are with respect to the nonconvex constraint set of the AC model, we first find the closest projection of the dispatch by solving the following optimization problem:
\begin{subequations} \label{eq:proj}
\begin{align}
\overline{\boldsymbol{p}}^{\rm g} \in \argmin_{\Xi^{\rm{AC}}} \quad & \sum_{k\in \mathcal{G}}\left(p_k^{\rm g} - \hat{p}_k^{\rm g} \right)^2 \label{eq:proj_a} \\
\textrm{s.t.\quad}   & \eqref{eq:1b}-\eqref{eq:1l},
\end{align}
\end{subequations}
where $\hat{p}_k^{\rm g}$ denotes the dispatch obtained from solving the (p)DC-OPF models, and $\overline{\boldsymbol{p}}^{\rm g}$ denotes the orthogonal projection onto the nonconvex constraint set $\eqref{eq:1b}$-$\eqref{eq:1l}$. Then, we define the distance to feasibility, which corresponds to the root mean square error between $\hat{\boldsymbol{p}}^{\rm g}$ and its orthogonal projection $\overline{\boldsymbol{p}}^{\rm g}$, as follows:
\[
\text{Distance to feasibility}\,(\textrm{p.u.}) = \sqrt{\frac{1}{\lvert \mathcal{G} \rvert} \sum_{k\in \mathcal{G}}\left(\overline{p}_k^{\rm g} - \hat{p}_k^{\rm g} \right)^2}.
\]
Figure~\ref{fig:Feasibility} depicts the distance to feasibility of the pDC- and DC-OPF models for the validation samples. Note that Problem \eqref{eq:proj} is solved using a local solver, which produces feasible solutions but not necessarily the closest feasible projection. The proposed model shows significant improvements in terms of dispatch feasibility compared to the DC-OPF model. In particular, for the IEEE 30-bus system, the improvement of the distance to feasibility is between 1 and 5 orders of magnitude with a maximum distance of less than $0.01$ p.u.

\subsection{Market Operations Analysis}
To assess the viability of the proposed model for market operations, we empirically study two elements: (1) LMPs and (2) market properties.

\begin{figure}[ht]
    \centering
    \includegraphics[width=\linewidth]{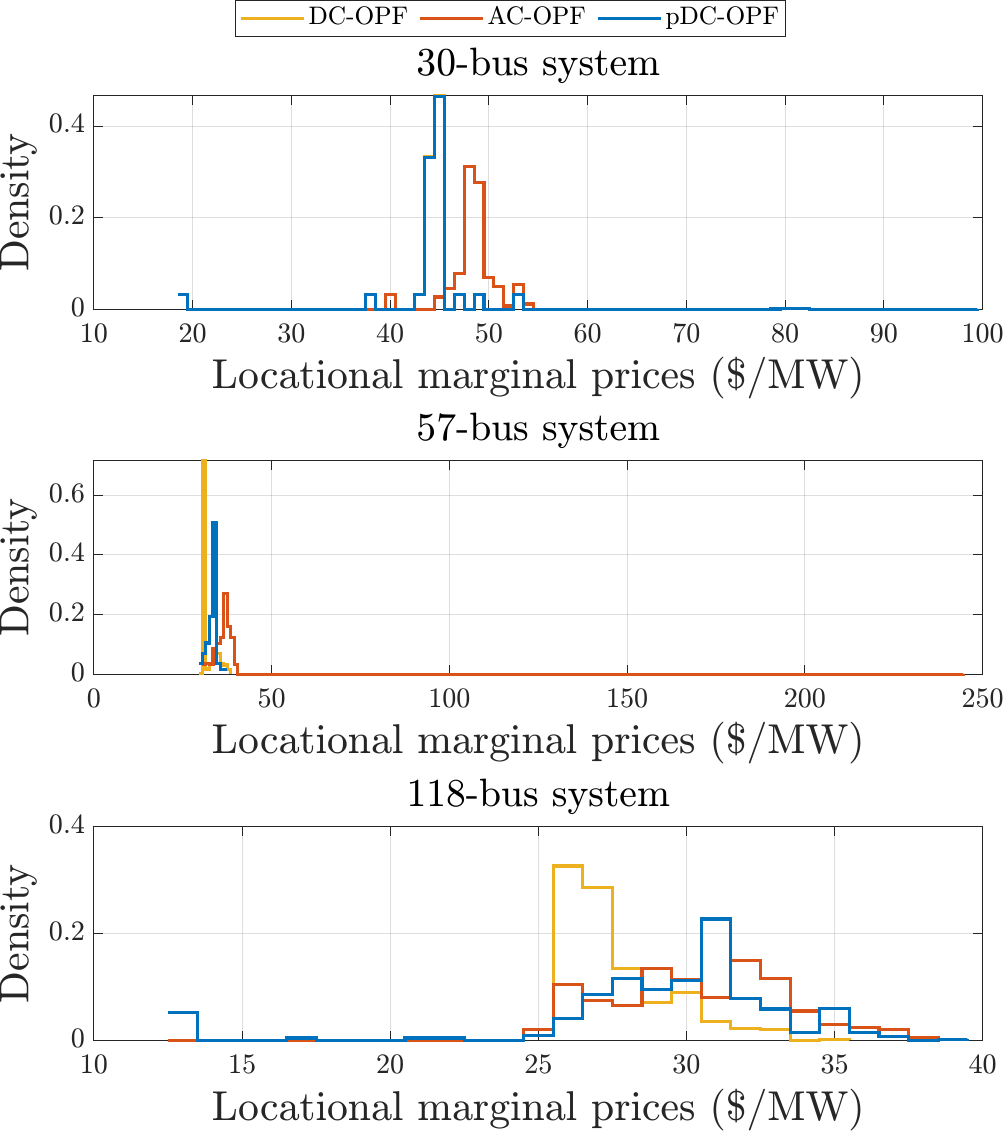}
    \caption{Locational marginal price distribution.}
    \label{fig:LMPS}
\end{figure}

Figure~\ref{fig:LMPS} depicts the distribution of the LMPs for the AC-, DC-, and pDC-OPF models. On average, the LMPs from the DC-OPF model underestimate the ones of the AC counterpart, which could result in significant financial implications, including distorted market signals for investments. That is, LMPs from the DC-OPF model might not necessarily reflect the system needs. This fact emphasizes the importance of using models that can accurately reflect the physics and constraints of the transmission system, whose implications go beyond the sole generation dispatch. 

For the 30-bus test system, both DC models render the same LMPs because the system has only two marginal generating units, i.e., both generators do not reach their upper or lower limits. The LMPs of the AC model of the 57-bus system show a significant tail corresponding to nodes where (i) voltage and reactive power limits are binding and (ii) lines that are congested due to considerable reactive power flows, which the (p)DC models cannot capture.

Additionally, we verify that the proposed model satisfies cost recovery and revenue adequacy properties in the online stage. Note that in the offline stage, the optimal scaling vector renders a dispatch and LMPs that satisfy both properties; however, due to the approximation error of the NN, such properties might not be satisfied at inference time. Table~\ref{tab:OkSamplesCR-RA} shows the fraction of samples that satisfy both properties when the optimal scaling vectors in the offline stage are determined enforcing or not such properties. The unconstrained model, which neglects constraints \eqref{eq:7b}-\eqref{eq:7c}, shows a higher fraction of samples not satisfying the two properties for the 57- and 118-bus systems, $4.1\%$ and $21\%$, respectively, showing the importance of considering such properties when generating the optimal scaling vectors.

\begin{table}[ht]
\centering
\caption{Verification of Market Properties}
\label{tab:OkSamplesCR-RA}
\begin{tabular}{c|c|ccc}
\hline
\multirow{2}{*}{\textbf{Model}} & \textbf{Market} & \multicolumn{3}{c}{\textbf{Fraction (\%)}} \\
                 & \textbf{Property} & \textbf{30-bus} & \textbf{57-bus} & \textbf{118-bus} \\
\hline
\multirow{2}{*}{\textbf{Original}} & RA  & 100   & 100   & 97.6    \\
                                      & CR & 100   & 100   & 100  \\
\hline
\multirow{2}{*}{\textbf{Unconst.}} & RA  & 100   & 95.9  & 79.0    \\
                                      & CR & 100   & 100   & 100  \\
\hline
\multicolumn{5}{l}{\footnotesize{RA: Revenue adequacy. CR: Cost Recovery.}}
\end{tabular}
\end{table}

\subsection{Topology Analysis}

In the offline stage, the AC-OPF solutions and the optimal scaling vectors are determined for a given network topology. However, the grid's topology often changes in practical settings due to planned equipment maintenance or unplanned failures/contingencies. We assess the proposed model's generalizability for network changes by testing its performance for N-1 network topologies, which are not represented in the dataset of the offline stage. In our experiments, we consider all N-1 network topologies that do not lead to islanded operation, and each topology is tested for the 1,000 demand samples.

\begin{figure}[t]
    \centering
    \includegraphics[width=\linewidth]{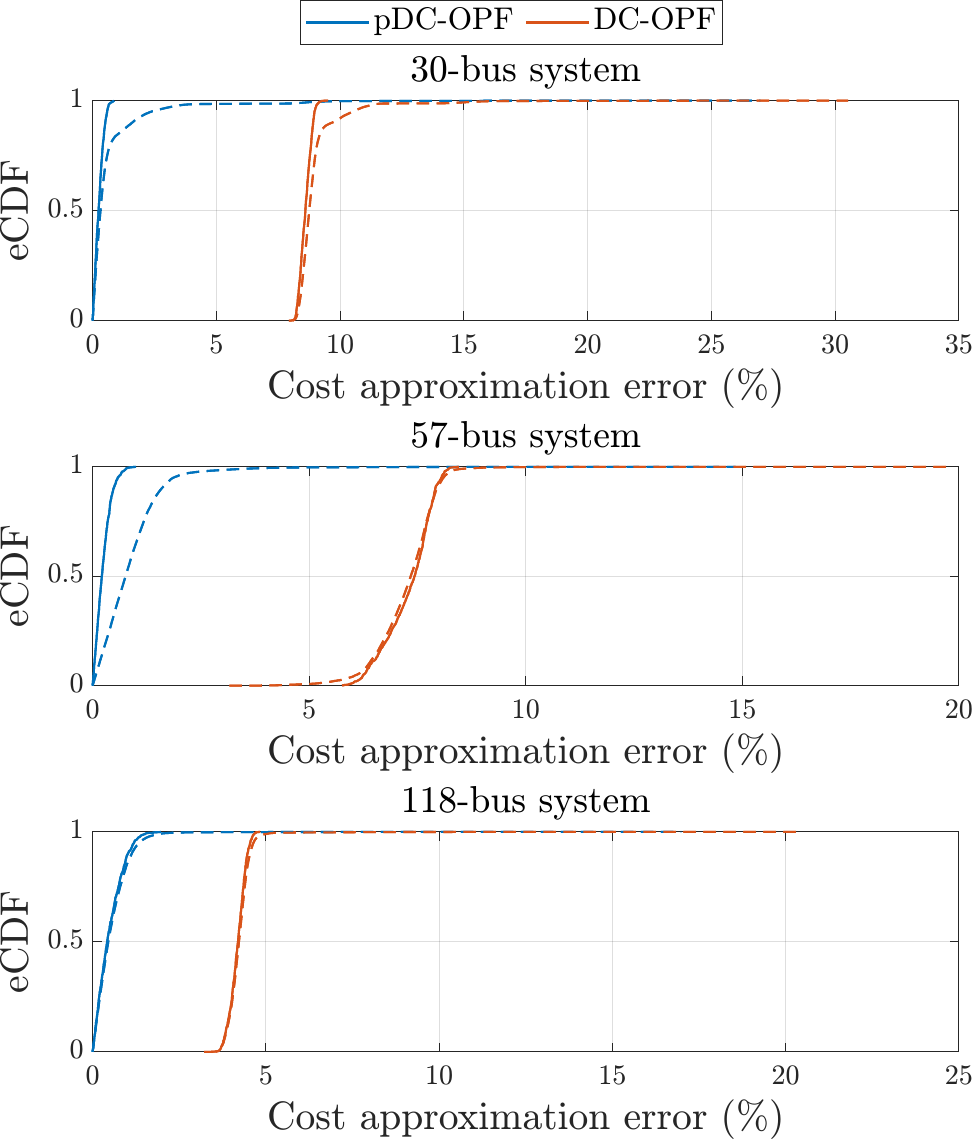}
    \caption{Cost approximation error empirical CDF for unrepresented N-1 topologies. Original topology (solid). N-1 topologies (dashed).}
    \label{fig:Contingency}
\end{figure}

Figure~\ref{fig:Contingency} depicts the empirical cumulative distribution function (eCDF) of the cost approximation error \eqref{eq:cost_error} of the original network topology and the N-1 network topologies for the pDC- and DC-OPF models. The cost approximation error of the proposed model with no contingencies for 98\% of the samples for the 30-, 57-, and 118-bus systems are less than $0.649\%$, $0.712\%$, and $1.382\%$, respectively, whereas for the N-1 topologies the errors are less than $3.602\%$, $3.171\%$, and $1.650\%$, respectively. Although the cost approximation errors for the N-1 topologies increase for the three systems, there is a clear benefit of the proposed model with respect to its DC counterpart. We note that the worsening of the errors for the N-1 topologies, with respect to the ones for the original topology, decreases as the network size increases, i.e., the proposed model's performance is less sensitive to the topological changes for larger systems.

\section{Concluding Remarks} \label{sec:conclusion}

This work presents a parametric optimization model to improve the dispatch accuracy of the DC-OPF model for market operations. The proposed model relies on a supervised learning framework to predict the optimal tunable parameters that reduce dispatch errors while satisfying market properties. Our numerical experiments show the benefits of the proposed model in terms of reduced cost errors, and dispatch accuracy and feasibility. Additionally, our results show that, on average, the locational marginal prices from the DC-OPF problem underestimate the ones of its AC counterpart, which could result in significant financial implications and highlights the importance of developing tractable models that provide a better representation of the power flow equations and the network constraints.

Our future work aims to improve the scalability of the proposed approach by (i) leveraging automatic differentiation tools and embedding the parameter determination problem into an end-to-end differentiable pipeline and (ii) approximately solving the bilevel problem using sequential convex programming \parencite{Constante2022} and linear reformulations derived from the KKT conditions of the lower level problem \parencite{Ruiz2009}.
Additionally, we plan to explore further parametric linear models that efficiently capture the impact of voltage and reactive power constraints on the LMPs and apply such models to other applications, such as unit commitment with security constraints and AC power flow equations.

\section{Acknowledgements}
This study was supported by startup funding from the Davidson School of Chemical Engineering and the College of Engineering at Purdue,  and the Office of Naval Research (grant no. N000142412641).

\printbibliography

\end{document}

%% file: hicss-packages.tex
\usepackage[letterpaper]{geometry}
\usepackage{hicss}
\usepackage{times}
\usepackage[none]{hyphenat}
\usepackage{url}
\usepackage{latexsym}
\usepackage{lineno}
\usepackage{indentfirst}
\usepackage{graphicx}
\graphicspath{{images/}}
\usepackage[
    style=apa,
  ]{biblatex}